\documentclass[11pt]{article}
\usepackage{amssymb}
\usepackage{amsmath}
\usepackage[dvips]{graphicx}

\input{tcilatex}
\begin{document}

\title{\textbf{ENHANCING THE ACCURACY OF THE TAYLOR POLYNOMIAL BY
DETERMINING THE REMAINDER TERM}}
\author{J.S.C Prentice\thanks{
\texttt{jpmsro@mathsophical.com}}}
\date{~}
\maketitle

\begin{abstract}
We determine the Lagrange function in Taylor polynomial approximation by
solving an appropriate initial-value problem. Hence, we determine the
remainder term which we then approximate by means of a natural cubic spline.
This results in a significant improvement in the quality of the Taylor
approximation. We observe improvements in the accuracy of the approximation
of many orders of magnitude, including a case when the independent variable $%
x$ lies beyond the relevant radius of convergence.

\textit{Key words: }Taylor, approximation, remainder, Lagrange function,
cubic spline
\end{abstract}

\pagestyle{myheadings} %Example
\markboth{Author1 Author2}{ENHANCING THE ACCURACY OF THE TAYLOR POLYNOMIAL} 
%<---SHORT TITLE

\section{Introduction}

In a recent paper \cite{Prentice 1}, we determined the remainder term in
Lagrange interpolation by solving suitable initial-value problem. In this
paper, we apply the same notion to Taylor approximation. We derive a
differential equation that will allow us to compute the reminder term (in
Lagrange form) of a first-order Taylor expansion. We then consider how a
cubic spline can be used to approximate the remainder term, leading to a
significant improvement in the accuracy of the approximation.

\section{Relevant Concepts}

Let $y\left( x\right) $ be a real-valued function, and assume that $y\left(
x\right) $ is as differentiable as is required in this paper (we will make
comments regarding this point when necessary). Taylor's theorem \cite{Spivak}
provides the following result:%
\begin{equation}
y\left( x\right) =y\left( x_{0}\right) +y^{\prime }\left( x_{0}\right)
\left( x-x_{0}\right) +\frac{y^{\prime \prime }\left( \xi _{x}\right) }{2}%
\left( x-x_{0}\right) ^{2}  \label{Taylor}
\end{equation}%
where $x_{0}<\xi _{x}<x.$ The third term on the RHS is the \textit{remainder
term}, presented here in \textit{Lagrange} form \cite{Apost}. There are
other representations of the remainder term, such as the Cauchy form \cite%
{Apost}%
\begin{equation*}
y^{\prime \prime }\left( \xi _{C}\right) \left( x-\xi _{C}\right) \left(
x-x_{0}\right) ,\text{ \ }x_{0}<\xi _{C}<x
\end{equation*}%
and the more general form%
\begin{equation}
y^{\prime \prime }\left( \xi _{C}\right) \left( x-\xi _{g}\right) \left(
x-x_{0}\right) \left( \frac{g\left( x\right) -g\left( x_{0}\right) }{%
g^{\prime }\left( \xi _{g}\right) }\right) ,\text{ \ }x_{0}<\xi _{g}<x
\label{general rem term}
\end{equation}%
of which the Lagrange and Cauchy forms are particular cases. In (\ref%
{general rem term}), the function $g$ is continuous on $\left[ x_{0},x\right]
$ and differentiable with a non-vanishing derivative on $\left(
x_{0},x\right) .$ However, it is the Lagrange form that will most easily
serve our purposes here.

If we define%
\begin{equation*}
T_{1}\left( x\right) \equiv y\left( x_{0}\right) +y^{\prime }\left(
x_{0}\right) \left( x-x_{0}\right)
\end{equation*}%
to be the so-called \textit{Taylor polynomial of first degree}, then we have%
\begin{equation}
y\left( x\right) -T_{1}\left( x\right) =\frac{y^{\prime \prime }\left( \xi
_{x}\right) }{2}\left( x-x_{0}\right) ^{2}.  \label{y(x)-T1(x)=...}
\end{equation}

\section{The Lagrange Function}

By differentiating (\ref{y(x)-T1(x)=...}) with respect to $x,$we find%
\begin{equation}
y^{\prime }\left( x\right) -T_{1}^{\prime }\left( x\right) =y^{\prime \prime
}\left( \xi _{x}\right) \left( x-x_{0}\right) +\frac{y^{\prime \prime \prime
}\left( \xi _{x}\right) \left( x-x_{0}\right) ^{2}}{2}\frac{d\xi _{x}}{dx}.
\label{dy-dyx0}
\end{equation}%
Since the LHS of this expression exists, we must assume that the RHS also
exists. Hence, we must assume that $\frac{d\xi _{x}}{dx}$ exists. This
implies that $\xi _{x}$ is a function of $x,$ and we will sometimes write $%
\xi \left( x\right) $ or simply $\xi $ in place of $\xi _{x}.$ We will refer
to $\xi \left( x\right) $ as the \textit{Lagrange function}. In (\ref%
{dy-dyx0}), we have used the notation%
\begin{eqnarray*}
y^{\prime \prime }\left( \xi \right) &\equiv &y^{\prime \prime }\left( \xi
\left( x\right) \right) \\
y^{\prime \prime \prime }\left( \xi \right) &\equiv &\frac{dy^{\prime \prime
}\left( \xi \left( x\right) \right) }{d\xi \left( x\right) }=\left. \frac{%
dy^{\prime \prime }\left( x\right) }{dx}\right\vert _{x=\xi }.
\end{eqnarray*}%
Rearranging (\ref{dy-dyx0}) and using $T_{1}^{\prime }\left( x\right)
=y^{\prime }\left( x_{0}\right) $ gives%
\begin{equation}
\frac{d\xi }{dx}=\frac{2\left( y^{\prime }\left( x\right) -y^{\prime }\left(
x_{0}\right) -y^{\prime \prime }\left( \xi \right) \left( x-x_{0}\right)
\right) }{y^{\prime \prime \prime }\left( \xi \right) \left( x-x_{0}\right)
^{2}}.  \label{dzeta/dx = ...}
\end{equation}

This is an initial-value problem that can, in principle, be solved to yield
the Lagrange function $\xi \left( x\right) $ for a suitable initial value.
Once $\xi \left( x\right) $ is known, the remainder term is easily computed.
Hence, if $T_{1}\left( x\right) $ is taken as an approximation to $y\left(
x\right) ,$ the resulting approximation error is known. We could also obtain
differential equations for $\xi _{C}$ and $\xi _{g}$ in a similar fashion,
but we believe these would be more complex than (\ref{dzeta/dx = ...}), and
so we work only with $\xi $ in this paper - an application of Occam's Razor,
we suppose.

There is a subtle point to be made: in (\ref{Taylor}), the remainder is
derived using the Mean Value Theorem \cite{boman}, hence the constraint $%
x_{0}<\xi _{x}<x.$ Our assumption above that $\xi _{x}=\xi \left( x\right) $
may be at odds with this constraint. The Lagrange form of the remainder term
does not explicitly require that $\xi _{x}$ be a smooth function of $x.$
This is not a concern for us. Our task here is simply to find values for $%
\xi $ that satisfy (\ref{Taylor}), by solving (\ref{dzeta/dx = ...}),
whether or not they satisfy the constraint. However, we will demonstrate
later how the solution(s) of (\ref{dzeta/dx = ...}) can be made consistent
with the constraint.

Lastly, note that our analysis up to this point requires that $y\left(
x\right) $ is three times differentiable, at least.

\section{Numerical Examples}

\subsection{First example}

For our first example, we consider%
\begin{equation*}
y\left( x\right) =e^{x/5}\sin x
\end{equation*}%
for $x\in \left[ 1,10\right] .$ So we have $x_{0}=1.$ To find a suitable
initial value for $\xi ,$ we choose a point very close to $x_{0}\ $(which we
denote $x_{z}$), say $x_{z}=1.0005$. We then solve 
\begin{equation}
y\left( x_{z}\right) =y\left( x_{0}\right) +y^{\prime }\left( x_{0}\right)
\left( x_{z}-x_{0}\right) +\frac{y^{\prime \prime }\left( \xi _{z}\right) }{2%
}\left( x_{z}-x_{0}\right) ^{2}  \label{Taylor xz}
\end{equation}%
numerically to find $\xi _{z}=1.000167.$ Next, we use the initial value $%
\left( x_{z},\xi _{z}\right) $ to solve (\ref{dzeta/dx = ...}) using a
seventh-order Runge-Kutta (RK7) method \cite{Butcher} to find $\xi \left(
x\right) .$ This then allows us to find the remainder%
\begin{equation*}
R_{\xi }\left( x\right) \equiv \frac{y^{\prime \prime }\left( \xi \left(
x\right) \right) }{2}\left( x-x_{0}\right) ^{2},
\end{equation*}%
which can be compared with the actual remainder%
\begin{equation*}
R_{act}\left( x\right) \equiv y\left( x\right) -T_{1}\left( x\right)
=y\left( x\right) -y\left( x_{0}\right) -y^{\prime }\left( x_{0}\right)
\left( x-x_{0}\right) .
\end{equation*}%
We can measure the quality of $R_{\xi }\left( x\right) $ through the device
of 
\begin{equation*}
\Delta R\left( x\right) \equiv R_{act}\left( x\right) -R_{\xi }\left(
x\right) .
\end{equation*}%
For this example, we find%
\begin{equation*}
\max_{\left[ 1,10\right] }\left\vert \Delta R\left( x\right) \right\vert
=5.4\times 10^{-13},
\end{equation*}%
suggesting that our calculation of $R_{\xi }\left( x\right) $ is very
accurate, indeed.

It transpires that (\ref{Taylor xz}) has another solution, leading to $\xi
_{z}=3.157781.$ We use this initial value in the same way as before,
eventually finding 
\begin{equation*}
\max_{\left[ 1,10\right] }\left\vert \Delta R\left( x\right) \right\vert
=3.0\times 10^{-13},
\end{equation*}%
which shows accuracy similar to the first case.

In Figure 1, we show $\xi \left( x\right) $ for both initial values. On this
plot we show the constant $x_{0}=1$ and the line $y=x.$ These are the bounds
on the constraint $x_{0}<\xi _{x}<x,x\in \left[ 1,10\right] .$ For $\left(
x_{z},\xi _{z}\right) =\left( 1.0005,1.000167\right) ,$ we see that $\xi
\left( x\right) $ lies between the two bounds up to $x\approx 5.4,$ and lies
below $x_{0}=1$ thereafter. By contrast, for $\left( x_{z},\xi _{z}\right)
=\left( 1.0005,1.000167\right) ,$ $\xi \left( x\right) >x$ up to $x\approx
2.8$ and thereafter lies between the bounds. Neither of these Lagrange
functions satisfies the constraint on the entire interval $\left[ 1,10\right]
.$ However, the Lagrange form of the remainder term, derived using the Mean
Value Theorem, does not require that $\xi _{x}$ must be smooth or even
continuous. It simply states that, for each $x\in \left[ 1,10\right] ,$
there does exist an $\xi _{x}$ such that $x_{0}<\xi _{x}<x$ and $\xi _{x}$
satisfies (\ref{Taylor xz}). If we choose values of $\xi \left( x\right) $
from the first Lagrange function for $x\in \left[ 1,4\right] ,$ and values
of $\xi \left( x\right) $ from the second Lagrange function for $x\in \left(
4,10\right] ,$ we necessarily satisfy Taylor's theorem for this example.

In Figure 2 we show $R_{\xi }\left( x\right) $ and $R_{act}\left( x\right) ,$
and in Figure 3 we show $\Delta R\left( x\right) .$

\subsection{Second example}

For our second example, we consider%
\begin{equation*}
y\left( x\right) =\ln \left( 1+x\right)
\end{equation*}%
for $x\in \left[ 0,10\right] .$ We have $x_{0}=0.$ Again, we choose $%
x_{z}=1.0005$. We then solve (\ref{Taylor xz}) to find $\xi _{z}=1.67\times
10^{-4}.$ Application of RK7 then gives the Lagrange function $\xi \left(
x\right) ,$ shown in Figure 4, along with the bounds of the constraint $%
x_{0}=0$ and the line $y=x.$ We see that $\xi \left( x\right) $ lies between
the bounds for all $x\in \left[ 0,10\right] .$ In Figure 5 we show $R_{\xi
}\left( x\right) $ and $R_{act}\left( x\right) ,$ and in Figure 6 we show $%
\Delta R\left( x\right) .$ We determine 
\begin{equation*}
\max_{\left[ 0,10\right] }\left\vert \Delta R\left( x\right) \right\vert
=1.3\times 10^{-13}.
\end{equation*}

\section{Polynomial Approximation of the Remainder}

Now that we have determined the Lagrange function $\xi \left( x\right) $
and, hence, the\ remainder term $R_{\xi }\left( x\right) $ for each example,
it seems reasonable to attempt to approximate the remainder term by means of
a polynomial. This polynomial can be added to the original Taylor polynomial 
$T_{1}\left( x\right) $ to yield a more accurate approximation than $%
T_{1}\left( x\right) .$ If we denote the polynomial approximation to $R_{\xi
}\left( x\right) $ by $P_{R}\left( x\right) ,$ we would then present%
\begin{equation*}
T_{1}\left( x\right) +P_{R}\left( x\right)
\end{equation*}%
as an approximation to $y\left( x\right) $ on the given interval.

We choose to use a \textit{natural} cubic spline \cite{B and F} to form $%
P_{R}\left( x\right) ,$ for several reasons: we have the RK nodes at our
disposal (we used $10000$ nodes in the RK computations); a cubic polynomial
combines with the factor $\left( x-x_{0}\right) ^{2}$ in $R_{\xi }\left(
x\right) $ to yield a polynomial of degree five, at most; and such splines
can be generated very efficiently on our computational platform \cite%
{platform}. Also, it is possible to estimate a bound on the accuracy of a
natural spline although, as will be seen, the bound is not tight.

Defining%
\begin{align*}
\Delta _{T}& \equiv \max_{I}\left\vert y\left( x\right) -T_{5}\left(
x\right) \right\vert \\
\Delta _{CS}& \equiv \max_{I}\left\vert y\left( x\right) -\left( T_{1}\left(
x\right) +P_{R}\left( x\right) \right) \right\vert
\end{align*}%
where $T_{5}\left( x\right) $ is the Taylor polynomial of fifth degree, and $%
I$ is the relevant interval of approximation, we show results in Table 1.%
\begin{equation*}
\end{equation*}

\begin{center}
\begin{tabular}{|c|c|c|c|c|}
\hline
$y\left( x\right) $ & $I$ & $\Delta _{T}$ & $\Delta _{CS}$ & $B_{U}$ \\ 
\hline
$e^{x/5}\sin x$ & $\left[ 1,10\right] $ & $5.8\times 10^{2}$ & $5.4\times
10^{-13}$ & $5.1\times 10^{-10}$ \\ \hline
$\ln \left( 1+x\right) $ & $\left[ 0,10\right] $ & $1.8\times 10^{4}$ & $%
1.4\times 10^{-13}$ & $8.6\times 10^{-9}$ \\ \hline
\end{tabular}

Table 1: Results for cubic spline and Taylor approximations.
\end{center}

\begin{equation*}
\end{equation*}%
In Table 1, $B_{U}$ is an upper bound on the error in the cubic spline
approximation, estimated using \cite{B and de}%
\begin{equation}
B_{U}=72h^{4}\max_{I}\left\vert y^{\left( 6\right) }\left( x\right)
\right\vert  \label{BU}
\end{equation}%
where $h$ is the uniform RK stepsize (see Appendix). The bounds are clearly
generous and, for each example, the actual error is less than the bound. It
is abundantly clear that $T_{1}\left( x\right) +P_{R}\left( x\right) $ is a
far better approximant than $T_{5}\left( x\right) .$ This is particularly
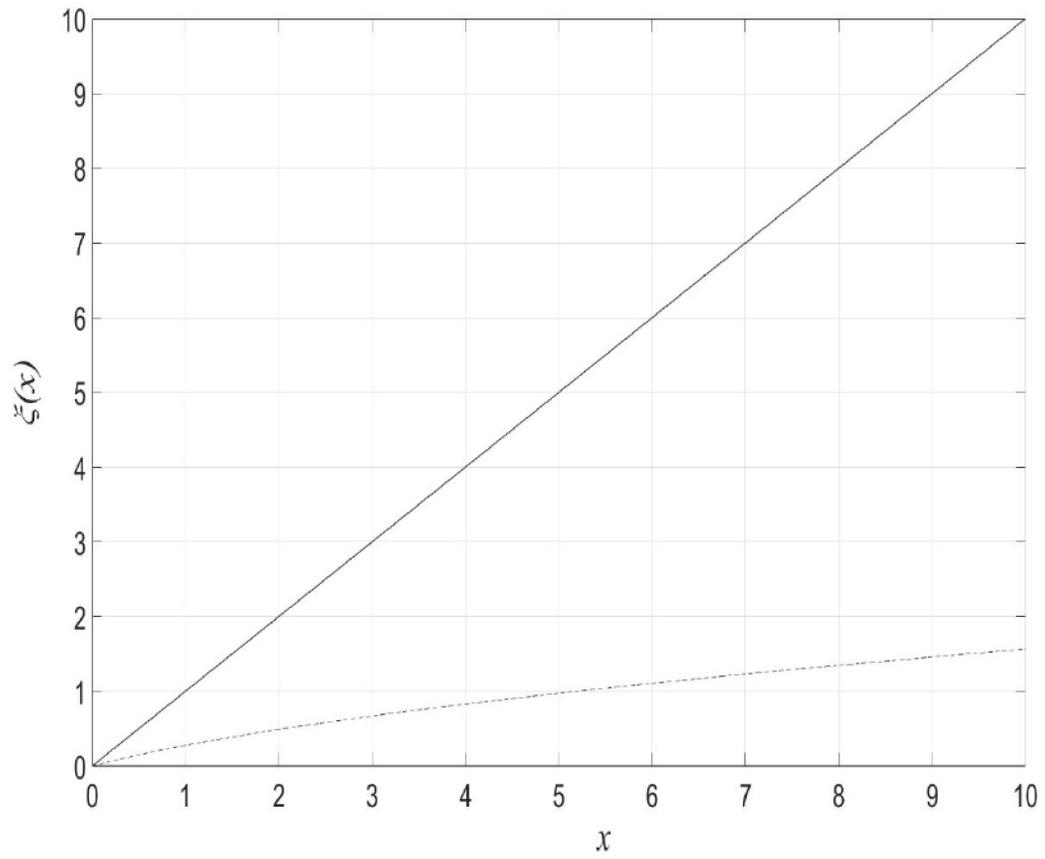
true for the second example, whose Taylor series has a radius of convergence
of $\left\vert x\right\vert <1$ (hence the very large value of $\Delta
_{T}). $

\section{Concluding Comments}

We have shown how the Lagrange function in Taylor polynomial approximation
can be determined by solving an appropriate initial-value problem. This
allows the remainder term to be determined. The remainder term can then be
approximated by means of a polynomial, and this can result in a significant
improvement in the quality of the Taylor approximation overall. We have
demonstrated this effect using a cubic spline, and we note improvements in
the accuracy of the approximation of many orders of magnitude, including the
case when the independent variable $x$ lies beyond the radius of
convergence. This speaks to the potential value of the idea presented here,
and in \cite{Prentice 1}.

This paper is intended as a demonstration of an idea and, as such, our
analysis has not been exhaustive. Further studies should consider the effect
of error control in the RK solution, and how this will affect the quality of
the polynomial approximation of the remainder term. It is also feasible to
approximate the remainder term by means of a least-squares fit, which might
be worth considering. Also, the multidimensional case must be investigated,
which would require the solution of a system of differential equations, but
all these aspects will be reserved for future work.

\medskip

\section{Appendix}

Here we derive the bound in (\ref{BU}). We refer to \cite{B and de}.
Following \cite{B and de} (see their eqn (27)), but using our own notation,
we have%
\begin{align*}
e\left( x\right) & \equiv y^{\prime \prime }\left( x\right) -P_{R}\left(
x\right) \\
\left\vert e^{\prime \prime \prime }\right\vert & =\left\vert y^{\left(
5\right) }\left( x\right) -P_{R}^{\prime \prime \prime }\left( x\right)
\right\vert \leqslant 3M\left( 1+M\right) ^{2}\max_{I}\left\vert \frac{%
d^{4}y^{\prime \prime }}{dx^{4}}\right\vert \max_{I}\left\vert
h_{i}\right\vert
\end{align*}%
where $h_{i}$ denotes the spacing between the RK nodes (which is not
necessarily uniform) and 
\begin{equation*}
M=\frac{\max_{I}\left\vert h_{i}\right\vert }{\min_{I}\left\vert
h_{i}\right\vert }.
\end{equation*}%
In our calculations, the RK nodes \textit{are} uniformly spaced (denoted $h$%
) so that $M=1.$ Hence,%
\begin{equation*}
\left\vert y^{\left( 5\right) }\left( x\right) -P_{R}^{\prime \prime \prime
}\left( x\right) \right\vert \leqslant 12\max_{I}\left\vert y^{\left(
6\right) }\right\vert h.
\end{equation*}%
Eqn (29) in \cite{B and de} gives the recursion%
\begin{align*}
\left\vert y^{\left( 4\right) }\left( x\right) -P_{R}^{\prime \prime }\left(
x\right) \right\vert & \leqslant \left( 3\right) \left\vert y^{\left(
5\right) }\left( x\right) -P_{R}^{\prime \prime \prime }\left( x\right)
\right\vert h \\
\left\vert y^{\left( 3\right) }\left( x\right) -P_{R}^{\prime }\left(
x\right) \right\vert & \leqslant \left( 2\right) \left\vert y^{\left(
4\right) }\left( x\right) -P_{R}^{\prime \prime }\left( x\right) \right\vert
h \\
\left\vert y^{\prime \prime }\left( x\right) -P_{R}\left( x\right)
\right\vert & \leqslant \left( 1\right) \left\vert y^{\left( 3\right)
}\left( x\right) -P_{R}^{\prime }\left( x\right) \right\vert h
\end{align*}%
so that%
\begin{align*}
\left\vert y^{\prime \prime }\left( x\right) -P_{R}\left( x\right)
\right\vert & \leqslant 6\left\vert y^{\left( 5\right) }\left( x\right)
-P_{R}^{\prime \prime \prime }\left( x\right) \right\vert h^{3} \\
& \leqslant 72\max_{I}\left\vert y^{\left( 6\right) }\right\vert h^{4}.
\end{align*}%
Note that the above analysis requires that $y\left( x\right) $ is six times
differentiable, at least.

\enddocument

###################################################

\end{document}